# Efficient modeling of chemotherapy regimens using mixed-integer linear programming

**Alain Billionnet**

Laboratoire CEDRIC, École Nationale Supérieure d'Informatique pour l'Industrie et l'Entreprise, 1, square de la Résistance, 91025, Évry cedex, France

email address : Alain.Billionnet@ensiie.fr

**Abstract**

In this article, we focus on determining a minimum-cost treatment program aimed at maintaining the size of a cancerous tumor at a level that allows the patient to live comfortably. At each predetermined point in a treatment horizon, the patient either receives drug treatment or does not. In the first case, the tumor shrinks and its size is multiplied by a constant factor lower than 1; in the second, it grows following an exponential or Gompertz growth law. We first demonstrate that a simple heuristic solution provides an optimal treatment program. We then show that the Gompertz function can be described, like the exponential function, by a simple recurrence relation that does not explicitly depend on time. Thanks to the characteristics of the logarithmic function, this property allows us to formulate the problem as a mixed-integer linear program. This result makes it possible to solve the problem very efficiently using one of the many solvers available to handle this type of program, and above all to consider and solve several extensions to the problem. In particular, we show how to determine which of the optimal equivalent solutions to the initial problem are the most relevant. We also show how to measure the effect of a marginal increase in the treatment budget on patient quality of life. Numerous computational experiments are presented to illustrate these issues.

## 1 Introduction

The problem studied in this article is inspired by the reference (Bräutigam, 2024). In this article, the author is interested in determining a chemotherapy treatment program for cancer that achieves a given medical objective while being cost-effective. Determining a treatment schedule involves deciding when the treatment should be administered to the patient. It is assumed that we are dealing with a tumor which, without treatment, grows and, on the contrary, regresses when treatment is administered. Both tumor growth and tumor regression in response to treatment are complex phenomena (Wodarz & Komarova, 2014). We assume that the models describing growth and regression are known deterministic models. This question of optimal treatment planning is part of cancer research, a disease which, despite major advances, remains a serious affliction. Readers can refer to (Bräutigam, 2024) and the references therein for many details and additions on the subject and also to (Shi et al., 2014), a comprehensive survey on optimization models concerning chemotherapy treatment planning in oncology. As discussed in (Brautigäm, 2024), the search for a good treatment schedule can be carried out heuristically. The result is a solution that may be quite satisfactory, but it is not clear whether it differs significantly from the optimal solution. In this article, we show that, under certain conditions, a heuristic solution that comes naturally to mind is optimal. An important issue in the search for appropriate solutions is that the problem under consideration generally admits several optimal solutions, i.e. of minimal cost, which are not at all equivalent from a medical point of view. We show in this article that the problem can be tackled quite effectively by mathematical programming, more specifically by mixed-integer linear programming. This makes it possible, on the one hand, to determine which of the optimal solutions is the most interesting according to certain criteria and, on the other hand, to consider several relevant extensions of the initial problem. It is worth noting that numerous articles have been published on how operations research can contribute to the fight against cancer. These articles address issues relating to prevention, diagnosis, staging and treatment (Saville et al., 2018).

## 2  Presentation of the problem

This work inspired by the reference (Bräutigam, 2024) concerns the scheduling of drug dose administration, over a certain treatment horizon, to cancer patients. This horizon is divided into $K$ periods of equal length (52 weeks, for example), at the start of which a drug dose may or may not be administered. It is assumed that we know how tumor size evolves over a period: if no treatment is administered at the start of the period, the size of the tumor increases over the period following a certain growth function; if a treatment is administered at the start of the period, the size of the tumor decreases following a certain decay function. More precisely, if we note $S_k$ the size of the tumor at the beginning of period $k$, we have $S_k = f(S_{k-1})$ if there is no treatment at the start of period $k$-1 and $S_k = g(S_{k-1})$ if there is treatment. It is assumed that $f(x) > x$ and $g(x) < x$ for all $x$ and that $f(x)$ and $g(x)$ are increasing functions. To simplify presentation, we make the assumption, as in (Bräutigam, 2024), that no treatment can be undertaken that would result in tumor size falling below a certain fixed threshold. We will comment on this assumption at the end of Section 4. The problem is to determine the periods at the start of which a dose of drug must be administered to maintain, over the entire horizon, the tumor size below a fixed threshold and thus ensure a certain quality of life for the patient while minimizing the overall cost generated by the treatment over the horizon considered. It is assumed here that the cost of one treatment is equal to *p*.

## 3 Under certain conditions, the heuristic solution is optimal

A natural heuristic solution is to administer treatment at the start of period *k* if and only if, without treatment, the tumor size at the start of period *k*+1 would exceed the tolerated size, $S_{Tol}$. In other words, if $\bar{S}$ verifies $f(\bar{S}) = S_{Tol}$, since the function *f* is increasing, the heuristic solution is to administer treatment at the start of period *k* if and only if $S_k > \bar{S}$. We'll show by recurrence on *k* that if $g(f(S)) \leq f(g(S))$ then the heuristic solution is optimal. Let $C_k[S]$ be the minimum treatment cost over periods *k*, *k*+1,..., *K* when the tumor size is equal to *S* at the beginning of period *k*. According to Bellman's optimality principle (Bellman, 1953),

$$C_k[S] = \min\{C_{k+1}[f(S)], p + C_{k+1}[g(S)]\} \quad (a)$$

with constraints $f(S) \leq S_{Tol}$ and $g(S) \geq S_{\min}$.

*The function $C_k[S]$ is increasing*

Since *f* and *g* are increasing functions, $C_K[S]$ is 0 if $S \leq \bar{S}$ and *p* if $S > \bar{S}$. $C_K(S)$ is therefore an increasing function of *S*. Assume that $C_{k+1}[S]$ is an increasing function and show that $C_k[S]$ is also increasing. According to (a), since *f*, *g* and $C_{k+1}$ are increasing functions of *S*, $C_k[S]$, which is equal to the minimum of two increasing functions, is also an increasing function.

*Optimality of the heuristic solution*

Recurrence assumption: At the start of period *k*+1, the optimal decision is to apply the heuristic, i.e. if $f(S) \leq S_{Tol}$ then $C_{k+1}[S] = C_{k+2}[f(S)]$. We deduce $C_{k+2}[f(S)] \leq p + C_{k+2}[g(S)]$. This is obvious at the start of the last period. Let's show that, at the start of period *k*, the optimal decision is also to apply the heuristic, i.e. if $f(S) \leq S_{Tol}$ then $C_k[S] = C_{k+1}[f(S)]$. To do this, we have to show that if $f(S) \leq S_{Tol}$ then $C_{k+1}[f(S)] \leq p + C_{k+1}[g(S)]$. According to the recurrence assumption, if $f(f(S)) \leq S_{Tol}$ (which also implies $f(g(S)) \leq S_{Tol}$ since *f* is increasing and $g(S) < f(S)$) then $C_{k+1}[f(S)] = C_{k+2}[f(f(S))] \leq p + C_{k+2}[g(f(S))]$ and $p + C_{k+1}[g(S)] = p + C_{k+2}[f(g(S))]$. The heuristic solution is therefore optimal if $g(f(S)) \leq f(g(S))$ since we've shown that $C_k$ is increasing. This property can be interpreted as follows: let's consider a tumor of a certain size at the start of period *k* and assume that the following 2 strategies are possible: 1) no treatment at the start of period *k* and treatment at the start of period *k*+1, 2) treatment at the start of period *k* and no treatment at the start of period *k*+1. If $g(f(S)) \leq f(g(S))$ then the tumor size at the start of period *k*+2 will be smaller with strategy 1 than with strategy 2. We'll see in Section 4 that this property is verified when the growth of the tumor size follows the exponential or Gompertz law, whereas its decrease consists in applying a multiplicative reduction factor.

## 4  General formulation of the problem by a mathematical program

Although the heuristic solution is an optimal solution in the case where $g(f(S)) \leq f(g(S))$, it is interesting to approach the problem by mathematical programming. Indeed, as we shall see later, many interesting extensions to the problem cannot be solved by a heuristic method, but can be solved by mathematical programming. Bräutigam (2024) proposes several formulations of the problem in the form of a mathematical program in the case where, on the one hand, tumor growth follows the exponential model (Collins, 1956) or the Gompertz model (Laird, 1964; Gerlee, 2013; Benzekry, 2014) and, on the other hand, the effect of treatment reduces tumor size by a constant factor. Bräutigam (2024) notes that solving these programs by the optimization solver CPLEX (IBM Inc, 2021) for the exponential model and by BARON (Sahinidis, 1996) for the Gompertz model requires considerable computing time. The use of these programs is therefore limited and Bräutigam (2024) points out that this *gives rise to the need for research to developing ways to improving computational performance*. In this section, we begin by expressing the tumor size at the start of period *k* as a general function of its size at the start of period *k*-1, based on the functions *f* and *g* from the previous section. This size at the start of period *k* naturally depends on the treatments administered at the start of periods 1,2,...,*k*-1. Then we will consider the case where $f(S_{k-1}) = \alpha(S_{k-1})^\beta$ and $g(S_{k-1}) = (1 - RF) f(S_{k-1})$ where $\alpha$ and $\beta$ are two positive coefficients and *RF* is the reduction factor associated with the treatment. As we will see in sections 6.1 and 6.2, the growth function *f* corresponds to the two cases considered by Bräutigam (2024): the exponential model when $\beta = 1$ and the Gompertz model when $0 < \beta < 1$. Thus, the evolution of tumor size between the beginning of period *k*-1 and the beginning of period *k* obeys a multiplicative recurrence with a possible power. The fact that this recurrence does not explicitly depend on *k* for either tumor growth or tumor decay will, as we shall see, make solving the problem much easier.

*Notation*

- **Data**

$H = \{1,...,K\}$ : horizon considered composed of $K$ periods of equal length (e.g. weeks);

$p$ : cost for a dose of treatment;

$S_{init}$: initial tumor size (at the beginning of period 1);

$S_{min}$ : minimum possible size of the tumor after treatment;

$S_{Tol}$ : maximum size of tumor tolerated;

$f$: function describing tumor growth without treatment, $S_{k+1} = f(S_k)$ ;

$g$: function describing tumor decay with treatment, $S_{k+1} = g(S_k)$ ;

$RF$: tumor size reduction coefficient after treatment, $S_{k+1} = (1-RF)f(S_k)$ ; in this case $g(S_k) = (1-RF)f(S_k)$.

- **Variables**

$X_k \in \{0,1\} \, (k \in H)$: Boolean variable that equals 1 if and only if a dose of treatment is administered at the beginning of period $k$;

$S_k \in \mathbb{R}^+, S_{min} \leq S_k \leq S_{Tol} \, (k \in H)$: variable representing tumor size at the beginning of period $k$, taking into account doses administered in previous periods; tumor size must always remain between two fixed extreme values, $S_{min}$ and $S_{Tol}$.

$LS_k \in \mathbb{R}^+, \log(S_{min}) \leq LS_k \leq \log(S_{Tol}) \, (k \in H)$: logarithm of $S_k$.

According to the above and using Boolean variables $X_k$ which by convention take the value 1 if and only if a treatment is administered at the beginning of period $k$, the tumor size at the start of period $k$, $S_k$, is equal to $(1-X_{k-1})f(S_{k-1}) + X_{k-1}g[f(S_{k-1})]$. Indeed, the tumor size at the start of period $k$ is equal to $g[f(S_{k-1})]$ if a dose of drug is administered at the start of period $k$-1, corresponding to $X_{k-1} = 1$, and to $f(S_{k-1})$ otherwise, corresponding to $X_{k-1} = 0$.

Now consider the case where $f(x) = \alpha x^\beta$ and $g(x) = (1-RF)x$. In this case, the previous expression for tumor size becomes $S_k = (1-RFX_{k-1})\alpha(S_{k-1})^\beta$. Indeed, the tumor size at the start of period $k$ is equal to $(1-RF)\alpha(S_{k-1})^\beta$ if a dose of drug is administered at the start of period $k$-1, corresponding to $X_{k-1} = 1$, and to $\alpha(S_{k-1})^\beta$ otherwise, corresponding to $X_{k-1} = 0$. Using the logarithmic function, we obtain the following equality: $\log(S_k) = \log(1 - RFX_{k-1}) + \log(\alpha) + \beta \log(S_{k-1})$. Examining the two possible values of the variable $X_{k-1}$, we see that $\log(1 - RFX_{k-1})$ is equal to $\log(1 - RF)$ if $X_{k-1} = 1$ and 0 otherwise. We therefore have $\log(1 - RFX_{k-1}) = X_{k-1} \log(1 - RF)$ and the previous equality is equivalent to the following equality $\log(S_k) = X_{k-1} \log(1 - RF) + \log(\alpha) + \beta \log(S_{k-1})$ or, by setting $LS_k = \log(S_k)$, to $LS_k = X_{k-1} \log(1 - RF) + \log(\alpha) + \beta LS_{k-1}$. The advantage of this last expression is its linear nature (as a function of the variables $LS_k$). The problem under consideration can therefore be formulated by the mixed-integer *linear* program P$_1$.

$$P_1: \begin{cases} \min \; C(X) = \sum_{k=1}^{K} p X_k \\ \text{s.t.} \begin{vmatrix} LS_1 = \log(S_{init}) & & (1.1) \\ LS_k = X_{k-1} \log(1 - RF) + \log(\alpha) + \beta LS_{k-1} & k \in H, k \geq 2 & (1.2) \\ LS_k \in \mathbb{R}^+, \log(S_{\min}) \leq LS_k \leq \log(S_{Tol}) & k \in H & (1.3) \\ X_k \in \{0,1\} & k \in H & (1.4) \end{vmatrix} \end{cases}$$

The economic function, $C(X)$, expresses the cost of treatment since, by convention, a dose of drug is administered at the start of period $k$ if and only if $X_k = 1$. Constraint 1.1 states that the logarithm of the tumor size at the start of period 1 is equal to $\log(S_{init})$. Linear constraint 1.2 expresses the logarithm of the tumor size at the start of period $k$ as a function of the logarithm of its size at period $k$-1, the variable $X_{k-1}$ and the coefficients $RF$, $\alpha$ and $\beta$. Constraints 1.3 and 1.4 specify the domain of variables $LS_k$ and $X_k$. In conclusion, the problem can be formulated by a linear mathematical program with real variables and Boolean variables in the case where 1) Tumor growth obeys the function $f(x) = \alpha x^\beta$ which is the case for the exponential or Gompertz model, 2) Tumor decay obeys the function $g(x) = (1 - RF)x$. This linear formulation will enable the problem to be solved efficiently using one of the many mixed-integer linear programming solvers available.

*Remark*

To simplify presentation, we have made the assumption, as in (Bräutigam, 2024), that no treatment can be undertaken that would reduce tumor size below a certain value. The following slightly different constraint could easily be taken into account: a treatment can always be undertaken at the beginning of each period, but the size reached by the tumor as a result of this treatment cannot fall below a certain value. To take account of this slight modification to the problem, it would suffice to replace the equality constraint (1.2) by the inequality constraint $LS_k = \max\{X_{k-1}\log(1-RF) + \log(\alpha) + \beta LS_{k-1}, \log(S_{\min})\}$. Some solvers, such as Gurobi (Gurobi Optimization, 2023), handle this type of constraint directly. If using a solver that doesn't, this constraint can be linearized beforehand in a classic way, using new Boolean variables.

## 5 Model strengthening

The fact that we can formulate the problem with mixed-integer *linear* programming in the case of the exponential model and the Gompertz model opens up the possibility of solving interesting extensions of the problem in these cases. We present several of these extensions below by way of example. We shall see that they are easy to formalize from the basic program $P_1$.

### 5.1 Combination therapy

Combination therapy, a treatment modality that combines two or more therapeutic agents, is a particularly interesting approach to cancer therapy. Indeed, this approach can have better efficacy than the mono-therapeutic approach due to synergistic or/and additive properties on the one hand, and reduce drug resistance on the other hand (Mokhtari et al., 2017). We will consider here the case where *n* types of tumor treatment are available and that one and only one of these treatments can potentially be administered at the start of each period. Let's assume that type *i* treatment (*i*=1,2,…,*n*) costs $p_i$ and corresponds to a reduction factor $RF^i$. Using the Boolean variable $X_k^i$ ($k=1,...,K; i=1,...,n$) which, by convention, is equal to 1 if and only if a treatment of type *i* is administered at the beginning of period *k*, we can express $S_k$ as a function of $S_{k-1}$ as $S_k = (1 - \sum_{i=1}^n RF^i X_{k-1}^i)\alpha(S_{k-1})^\beta$ as long as we impose the constraint $\sum_{i=1}^n X_{k-1}^i \leq 1$. Using the logarithmic function and noting the logarithm of $S_k$ by $LS_k$ as before, we obtain the following equality:

$$LS_k = \log(1 - \sum_{i=1}^n RF^i X_{k-1}^i) + \log(\alpha) + \beta LS_{k-1}.$$

Since $\sum_{i=1}^n X_{k-1}^i \leq 1$, $\log(1 - \sum_{i=1}^n RF^i X_{k-1}^i) = \sum_{i=1}^n \log(1-RF^i) X_{k-1}^i$. Finally, the problem can be formulated by the program $P_2$.

$$P_2: \begin{cases} \min C'(X) = \sum_{k=1}^{K}\sum_{i=1}^{n} p_i X_k^i \\ \text{s.t.} \begin{cases} LS_1 = \log(S_{init}) & & (2.1) \\ LS_k = \sum_{i=1}^{n}\log(1-RF^i)X_{k-1}^i + \log(\alpha) + \beta LS_{k-1} & k \in H, k \geq 2 & (2.2) \\ \sum_{i=1}^{n} X_{k-1}^i \leq 1 & k \in H & (2.3) \\ LS_k \in \mathbb{R}^+, \log(S_{\min}) \leq LS_k \leq \log(S_{Tol}) & k \in H & (2.4) \\ X_k^i \in \{0,1\} & k \in H, i = 1,...,n & (2.5) \end{cases} \end{cases}$$

## 5.2 Spacing between treatments

Time interval to treatment is an important question (Yoo et al., 2017). It can be interesting, say, to have solutions that ensure a minimum spacing between two treatments. For example, $\delta$ treatment-free periods can be imposed between two treatment periods. Thus, if a treatment takes place at the start of period $k$, there can be no treatment at the start of periods $k+1$, $k+2$, $k+\delta$. Of course, respecting these treatment intervals can increase the overall cost of treatment. Under this constraint, the problem can be formulated as program $P_3$.

$$P_3: \begin{cases} \min C(X) \\ \text{s.t.} \begin{cases} LS_1 = \log(S_{init}) & & (3.1) \\ LS_k = X_{k-1}\log(1-RF) + \log(\alpha) + \beta LS_{k-1} & k \in H, k \geq 2 & (3.2) \\ X_{k+1} + X_{k+2} + ... + X_{k+\delta} \leq \delta(1-X_k) & k = 1,...,K-\delta & (3.3) \\ LS_k \in \mathbb{R}^+, \log(S_{\min}) \leq LS_k \leq \log(S_{Tol}) & k \in H & (3.4) \\ X_k \in \{0,1\} & k \in H & (3.5) \end{cases} \end{cases}$$

Let's look at constraint (3.3). If a treatment is administered at the beginning of period $k$, then $X_k = 1$ and the second member of the constraint equals 0. This implies $X_{k+1} = X_{k+2} = ... = X_{k+\delta} = 0$. Otherwise, the second member of the constraint equals $\delta$ and the constraint becomes inactive.

Note that the formulation of the problem by the mathematical program $P_1$ also makes it easy to take into account the fact that certain treatment dates may be imposed a priori, for example at the beginning of the 10th week and at the beginning of the 20th. To do this, simply add the constraints $X_k = 1$ to the program for the relevant values of $k$.

## 5.3 Minimizing maximum tumor size

It may be of interest to determine a treatment schedule that minimizes the maximum tumor size over the time horizon under consideration, while remaining below a certain cost.

This question can be associated with each of the 3 preceding problems, formulated as $P_1$, $P_2$ and $P_3$, which themselves aim to minimize cost while ensuring that the tumor does not grow beyond a certain size. For $P_1$, $P_2$ and $P_3$ this problem can be formulated by replacing in these programs the economic function to be minimized by the non-negative real variable $LS$ and adding the constraint $LS_k \leq LS$ ($k \in H$). For $P_1$ and $P_3$, the constraint $\sum_{k \in H} pX_k \leq C$ must also be added, and for $P_2$ the constraint $\sum_{k=1}^{K} \sum_{i=1}^{n} p_i X_k^i \leq C$ where $C$ denotes the cost not to be exceeded. $\pi_1(C)$, $\pi_2(C)$ and $\pi_3(C)$ are the programs obtained.

$$\pi_1(C): \begin{cases} \min LS \\ \text{s.t.} \begin{vmatrix} (1.1)\ (1.2) \\ (1.3)\ (1.4) \\ \sum_{k=1}^{K} pX_k \leq C \\ LS \geq LS_k\ (k \in H) \\ LS \in \mathbb{R}^+ \end{vmatrix} \end{cases} \quad \pi_2(C): \begin{cases} \min LS \\ \text{s.t.} \begin{vmatrix} (2.1)\ (2.2)\ (2.3) \\ (2.4)\ (2.5) \\ \sum_{k=1}^{K} \sum_{i=1}^{n} p_i X_k^i \leq C \\ LS \geq LS_k\ (k \in H) \\ LS \in \mathbb{R}^+ \end{vmatrix} \end{cases} \quad \pi_3(C): \begin{cases} \min LS \\ \text{s.t.} \begin{vmatrix} (3.1)\ (3.2)\ (3.3) \\ (3.4)\ (3.5) \\ \sum_{k=1}^{K} pX_k \leq C \\ LS \geq LS_k\ (k \in H) \\ LS \in \mathbb{R}^+ \end{vmatrix} \end{cases}$$

In an optimal solution of $\pi_i(C)$, $i=1,2,3$, the maximum size reached by the tumor during horizon $H$ is equal to $e^{LS*}$ where $LS*$ designates the optimum value of $\pi_i(C)$, i.e. the optimum value of variable $LS$.

In addition, since the program $P_1$ (resp. $P_2$, $P_3$) generally has several optimal solutions, it may be of particular interest to determine from among these optimal solutions, the one that minimizes the maximum size taken by the tumor over the horizon under consideration. The program $\pi_1(C)$ (resp. $\pi_2(C), \pi_3(C)$) can be used to determine this particular optimal solution. To do this, simply solve $\pi_1(C)$ (resp. $\pi_2(C), \pi_3(C)$), giving $C$ the optimal value of $P_1$ (resp. $P_2$, $P_3$).

## 6 Experiments

In this section, we present experiments to solve the programs $P_i$ and $\pi_i(C)$, $i=1,2,3$. Section 6.1 deals with the exponential tumor growth model and section 6.2 with the Gompertz model. In these two sections, the effect of treatment is to reduce tumor size by a factor of $RF$. All the experiments presented in this article were carried out using the mathematical programming language AMPL (Fourer and Gay, 1993), version 2022-10-13, to model the different programs and Gurobi (Gurobi Optimization, 2023), version 10.0.1, a commercially available mathematical program solver based on the most efficient algorithms available today, to solve them. Gurobi can solve linear and non-linear programs involving both real and integer variables. The experiments have been performed on a PC with an Intel Core i7 1.90 GHz processor with 16 Go RAM.

### 6.1 Exponential model

In this section, we consider the case where tumor growth follows the exponential distribution given by: $\varphi(t) = \varphi_0 e^{bt}$, where $\varphi(t)$ is the value of the function $\varphi$ at time $t$, $\varphi_0$ is its value at time $t=0$, $b$ is a positive coefficient expressing the growth rate and $e$ is the base of the natural logarithm. It is generally accepted that in the early stages, tumor growth can be approximated by an exponential function (Retsky et al. 1990). It is assumed here and in all that follows that the unit of time is the period. It is known (and easy to demonstrate) that, in the exponential case, if the size of the tumor at the start of period $k$ is equal to $S_k$ then, without treatment, its size will be equal to $\alpha S_k$ at the start of period $k+1$ with $\alpha = \log(b)$. We are clearly within the framework considered, $S_{k+1} = \alpha (S_k)^\beta$, by posing $\alpha = \log(b)$ and $\beta = 1$. It's easy to check that in this case, based on the results of Section 3, the heuristic solution is optimal. This is because $g(f(S_k)) = f(g(S_k)) = \alpha^2 (1 - RF) S_k$.

Table 2 presents numerical experiments on solving $P_1$ (resp. $P_2$, $P_3$) and $\pi_1(C)$ (resp. $\pi_2(C)$, $\pi_3(C)$), giving $C$ the optimal value of $P_1$ (resp. $P_2$, $P_3$) for the latter three programs. The values of the various parameters used in the formulation of these programs, taken from (Bräutigam, 2024) with the exception of $RF_1$, $RF_2$, $p_1$ and $p_2$, are presented in Table 1. The horizon considered is composed of 52 periods.

**Table 1**
Parameter values for the treatment planning problem in the case of exponential tumor growth (Bräutigam, 2024)

| Tumor growth function parameters | $S_{init}$ | $S_{Tol}$ | $S_{min}$ | $RF$ | $RF_1$ | $RF_2$ | $p$ | $p_1$ | $p_2$ | $\delta$ |
|---|---|---|---|---|---|---|---|---|---|---|
| $\varphi_0 = 50$, $b = e^{1.5} \Rightarrow \alpha = \log(b) = 1.5$ | 50 | 500 | 10 | 0.6 | 0.6 | 0.7 | 10 | 10 | 13 | 1 |

**Comments on Table 2**

*First column*: Solving each of the 3 programs is instantaneous. It required less than one second of computation. The resolution of $P_1$ provides a treatment schedule with a cost of 210. With this schedule, the maximum size reached by the tumor during the horizon considered is equal to 496.46 and is therefore very close to the tolerated size ($S_{Tol}$=500). We may ask whether P$_1$ admits other optimal solutions and in particular solutions for which the maximum size reached by the tumor is smaller. Solving $\pi_1(C)$ with $C$=210 answers this question by providing a treatment schedule, still costing 210, but with a maximum tumor size equal to 315.48. This solution seems much more interesting than the first (at equal cost) in terms of the patient's quality of life, since it enables tumor size to be maintained at a much lower value. There is also the question of whether an additional treatment dose significantly influences the maximum tumor size reached. The resolution of $\pi_1(C)$ with $C$=220 answers this question by providing a treatment schedule that keeps the tumor size below a much smaller value (126.19). An additional dose of treatment is therefore very effective in this case, if the aim is to keep the tumor size as small as possible. These 3 solutions are shown in Figure 1. We can also see that solving P$_1$ provides fairly regular treatment frequencies, since between two treatment periods there is always at least one without treatment and at most two. This is not the case for the optimal solution of $\pi_1(C)$ with $C$=210 and also with $C$=220. In both cases, there are at least two consecutive periods with treatment and 4 or 6 consecutive periods without treatment between two treatment periods. Treatment scheduling is therefore much more irregular. Programs P$_3$ and $\pi_3(C)$ take this aspect of the problem into account.

*Second column*: The resolution of P$_2$ and $\pi_2(C)$ - with 2 possible types of treatment - is instantaneous. With the values chosen for the cost of treatments, the optimal treatment schedule includes, in the 3 cases examined (P$_2$, $\pi_2(204)$, and $\pi_2(217)$), a comparable number of type 1 and type 2 treatments. On the other hand, treatments are not administered on a regular basis. In all 3 cases, there are at least two consecutive periods with treatment and 4 or 6 consecutive periods without treatment between two treatment periods. The resolution of P$_2$ provides a solution of cost 204 with a maximum tumor size equal to 493.50. Resolving $\pi_2(204)$ does not reduce this size. On the other hand, if we accept an additional treatment, either type 1 or type 2, the maximum tumor size over the treatment horizon decreases considerably, from 493.50 to 148.05.

*Third column*: Resolution of P$_3$ and $\pi_3(C)$ is instantaneous. The cost of the optimal solution of P$_3$ and the maximum tumor size reached in this solution are identical to those of the solution of P$_1$. And so, as expected, imposing the additional constraint (compared with P$_1$) of having a treatment-free period between two treatment periods does not increase the value of the optimal solution or the maximum value of the tumor size. On the other hand, the addition of this constraint leads to a different solution to that obtained with $P_1$, insofar as there is greater irregularity in the administration of treatments (8 periods without treatment between two treatments). The solutions of $\pi_3(210)$ and $\pi_3(220)$ are identical to those of $\pi_1(210)$ and $\pi_1(220)$ in terms of maximum tumor size. On the other hand, it provides much greater regularity in treatment administration.

**Table 2**
Experimental results for solving the programs $P_1$, $P_2$, $P_3$, $\pi_1(C)$, $\pi_2(C)$ and $\pi_3(C)$ in the case of exponential tumor growth

| $P_1$ | $P_2$ | $P_3$ |
|---|---|---|
| Minimum cost: 210<br>Maximum tumor size: 496.46<br>Min. spacing between 2 treatments: 1<br>Max. spacing between 2 treatments: 2<br>Computation time: <1s.<br>Figure: 1 | Minimum cost: 204<br>Number of type 1 treatments: 10<br>Number of type 2 treatments: 8<br>Maximum tumor size: 493.50<br>Min. spacing between 2 treatments: 0<br>Max. spacing between 2 treatments: 5<br>Computation time: <1s. | Minimum cost: 210<br>Maximum tumor size: 496.46<br>Min. spacing between 2 treatments: 1<br>Max. spacing between 2 treatments: 8<br>Computation time: <1s. |
| $\pi_1(210)$ | $\pi_2(204)$ | $\pi_3(210)$ |
| Maximum tumor size: 315.48<br>Min. spacing between 2 treatments: 0<br>Max. spacing between 2 treatments: 6<br>Computation time: <1s.<br>Figure: 1 | Number of type 1 treatments: 10<br>Number of type 2 treatments: 8<br>Maximum tumor size: 493.50<br>Min. spacing between 2 treatments: 0<br>Max. spacing between 2 treatments: 5<br>Computation time: <1s. | Maximum tumor size: 315.48<br>Min. spacing between 2 treatments: 1<br>Max. spacing between 2 treatments: 3<br>Computation time: <1s. |
| $\pi_1(220)$ | $\pi_2(217)$ | $\pi_3(220)$ |
| Maximum tumor size: 126.19<br>Min. spacing between 2 treatments: 0<br>Max. spacing between 2 treatments: 4<br>Computation time: <1s.<br>Figure: 1 | Number of type 1 treatments: 10<br>Number of type 2 treatments: 9<br>Maximum tumor size: 148.05<br>Min. spacing between 2 treatments: 0<br>Max. spacing between 2 treatments: 6<br>Computation time: <1s. | Maximum tumor size: 126.19<br>Min. spacing between 2 treatments: 1<br>Max. spacing between 2 treatments: 2<br>Computation time: <1s. |

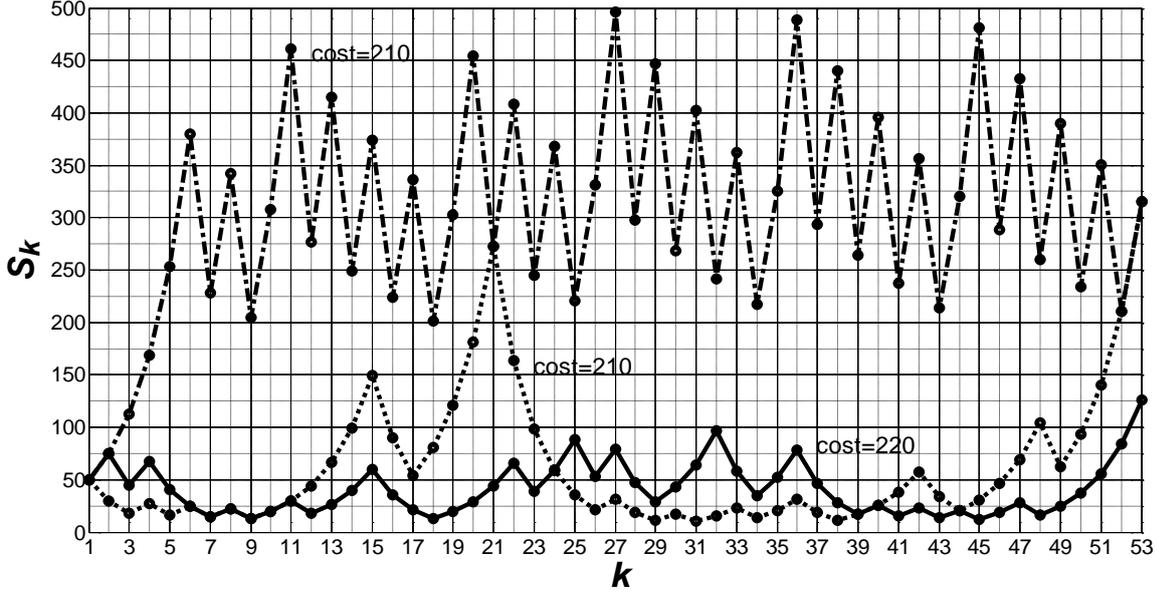

**Fig. 1** Tumor size at the beginning of each period for the optimal solution of $P_1$ (cost 210, dash-dotted line), for the optimal solution of $\pi_1(210)$ (dotted line) and for the optimal solution of $\pi_1(220)$ (solid line). In the solution of $P_1$, the maximum tumor size is equal to $S_{27}=496.46$, in the solution of $\pi_1(210)$ it is equal to $S_{53}=315.48$ and in the solution of $\pi_1(220)$ it is equal to $S_{53}=126.19$.

### 6.2 Gompertz model

In this section, we consider the case where tumor growth follows the Gompertz function. This model is frequently used to model tumor growth in biology and medicine as a function of time (Benzekry, 2014; Pitchaimani, 2014). It describes rapid initial tumor growth, which gradually slows as it approaches an asymptotic size. There are several equivalent expressions. We have adopted the one used by Bräutigam (2024): $\phi(t) = \phi_0\, e^{\frac{a}{b}(1-e^{-bt})}$ where $\phi(t)$ represents the value of the function at time $t$, $\phi_0$ its value at $t=0$ and $a$ and $b$ are two positive constants that control the growth rate and the shape of the associated curve.

Let us express $S_k$ as a function of $S_{k-1}$: $S_k / S_{k-1} = \phi_0\, e^{\frac{a}{b}(1-e^{-bk})} / \phi_0\, e^{\frac{a}{b}(1-e^{-b(k-1)})}$ $= e^{\frac{a}{b}(e^{-b(k-1)})(1-e^{-b})}$. Since $S_{k-1} = \phi_0\, e^{\frac{a}{b}} / e^{\frac{a}{b} e^{-b(k-1)}}$, we get $S_k = S_{k-1}(\phi_0\, e^{\frac{a}{b}} / S_{k-1})^{(1-e^{-b})} = (S_{k-1})^{e^{-b}}(\phi_0\, e^{\frac{a}{b}})^{(1-e^{-b})}$. We are therefore in the case where $f(S_{k-1})=\alpha(S_{k-1})^\beta$ by posing $\alpha = (\phi_0\, e^{\frac{a}{b}})^{(1-e^{-b})}$ and $\beta = e^{-b}$. Let's check that in this case the heuristic solution is optimal.

$$g(f(S_k)) = (1-RF)\alpha[\alpha(S_k)^{\beta}]^{\beta} = (1-RF)\alpha^{\beta+1}(S_k)^{\beta^2}, \qquad f(g(S_k)) = \alpha[(1-RF)\alpha(S_k)^{\beta}]^{\beta}$$

$= (1-RF)^{\beta}\alpha^{\beta+1}(Sk)^{\beta^2}$. We therefore have $g(f(S_k)) \leq f(g(S_k))$ if $1-RF \leq (1-RF)^{\beta}$, i.e. if $\log(1-RF) \leq \beta\log(1-RF)$. This last inequality is verified if $\beta \leq 1$, which is indeed the case since $\beta = e^{-b}$.

Table 4 presents numerical experiments on solving $P_1$ (resp. $P_2$, $P_3$) and $\pi_1(C)$ (resp. $\pi_2(C)$, $\pi_3(C)$) giving $C$ the optimal value of $P_1$ (resp. $P_2$, $P_3$) for the latter three programs. The values of the various parameters used in the formulation of these programs, taken from (Bräutigam, 2024) with the exception of $RF_1$, $RF_2$, $p_1$ and $p_2$, are presented in Table 3. The horizon considered is composed of 52 periods.

**Table 3**
Parameter values for Gompertz tumor growth experiments (Bräutigam, 2024)

| Tumor growth function parameters | $S_{init}$ | $S_{Tol}$ | $S_{min}$ | $RF$ | $RF_1$ | $RF_2$ | $p$ | $p_1$ | $p_2$ |
|---|---|---|---|---|---|---|---|---|---|
| $\phi_0 = 50$  $a=0.72$  $b=0.18$ $\Rightarrow \alpha = (g_0 e^{\frac{a}{b}})^{(1-e^{-b})} \approx 3.68, \beta = e^{-b} \approx 0.84$ | 150 | 500 | 60 | 0.6 | 0.6 | 0.7 | 10 | 10 | 13 |

**Table 4**
Experimental results for solving the programs $P_1$, $P_2$, $P_3$, $\pi_1(C)$, $\pi_2(C)$ and $\pi_3(C)$ in the case of tumor growth following Gompertz's law

| $P_1$ | $P_2$ | $P_3$ |
|---|---|---|
| Minimum cost: 200<br>Maximum tumor size: 499.09<br>Min. spacing between 2 treatments: 1<br>Max. spacing between 2 treatments: 2<br>Computation time: <1s.<br>Figure: 2 | Minimum cost: 191<br>Number of type 1 treatments: 10<br>Number of type 2 treatments: 7<br>Maximum tumor size: 499.09<br>Min. spacing between 2 treatments: 1<br>Max. spacing between 2 treatments: 3<br>Computation time: 29s. | Minimum cost: 200<br>Maximum tumor size: 498.87<br>Min. spacing between 2 treatments: 1<br>Max. spacing between 2 treatments: 3<br>Computation time: <1s. |
| $\pi_1(200)$ | $\pi_2(191)$ | $\pi_3(200)$ |
| Maximum tumor size: 438.42<br>Min. spacing between 2 treatments: 1<br>Max. spacing between 2 treatments: 2<br>Computation time: <1s.<br>Figure: 2 | Number of type 1 treatments: 10<br>Number of type 2 treatments: 7<br>Maximum tumor size: 493.41<br>Min. spacing between 2 treatments: 1<br>Max. spacing between 2 treatments: 2<br>Computation time: 36s. | Maximum tumor size: 438.42<br>Min. spacing between 2 treatments: 1<br>Max. spacing between 2 treatments: 2<br>Computation time: <1s. |
| $\pi_1(210)$ | $\pi_2(204)$ | $\pi_3(210)$ |
| Maximum tumor size: 418.02<br>Min. spacing between 2 treatments: 1<br>Max. spacing between 2 treatments: 2<br>Computation time: 3.2s.<br>Figure: 2 | Number of type 1 treatments: 19<br>Number of type 2 treatments: 1<br>Maximum tumor size: 435.78<br>Min. spacing between 2 treatments: 1<br>Max. spacing between 2 treatments: 2<br>Computation time: 70s. | Maximum tumor size: 418.02<br>Min. spacing between 2 treatments: 1<br>Max. spacing between 2 treatments: 2<br>Computation time: <1s. |

**Comments on Table 4**

*First column*: All 3 programs were solved very fast. Solving $P_1$ and $\pi_1(C)$ with $C=200$ took less than a second, while solving $\pi_1(C)$ with $C=210$ took 3.2 seconds. The resolution of $P_1$ provides a treatment schedule with a cost of 200. With this schedule, the maximum size reached by the tumor during the horizon considered is equal to 499.09 and is therefore roughly equal to the tolerated size ($S_{Tol}=500$). We may ask whether $P_1$ admits other optimal solutions and in particular solutions for which the maximum size reached by the tumor is smaller. Solving $\pi_1(C)$ with $C=200$ answers this question by providing a treatment schedule, still costing 200, but with a maximum tumor size equal to 438.42. This solution is slightly more attractive than the first (at equal cost) in terms of patient quality of life, since it enables tumor size to be maintained at a slightly lower value. There is also the question of whether an additional treatment dose significantly influences the maximum tumor size reached. The resolution of $\pi_1(C)$ with C=210 answers this question by providing a treatment regimen that keeps tumor size below a slightly smaller value (418.02). The value of an additional treatment dose is therefore limited in this case, if the aim is to keep tumor size as small as possible. These 3 solutions are shown in Figure 2. We can also see that solving $P_1$, $\pi_1(200)$ and $\pi_1(210)$ provides fairly regular treatment frequencies, since between two treatment periods there is always at least one without treatment and at most two.

*Second column*: Resolution of each of the 3 programs - with 2 possible types of treatment - is no longer instantaneous, but is still very fast. With the values chosen for treatment costs, the optimal treatment plan includes, in the 2 cases $P_2$ and $\pi_2(191)$, a relatively comparable number of type 1 and type 2 treatments. However, this is not the case for $\pi_2(204)$, since of the 20 treatments administered, only one is type 2. On the other hand, treatments are administered on a regular basis since between two treatment periods there is always at least one without treatment and at most two or three. The resolution of $P_2$ provides a solution of cost 191 with a maximum tumor size equal to 499.09. The resolution of $\pi_2(191)$ barely reduces this size. On the other hand, if we accept an additional treatment, either type 1 or type 2, the maximum tumor size over the treatment horizon decreases slightly, from 493.41 to 435.78.

*Third column*: Resolution of each of the 3 programs is instantaneous. We obtain exactly the same results as in column 1, which is not surprising since the 3 solutions in column 1 already respect the spacing constraint between two treatments imposed in programs $P_3$, $\pi_3(200)$ and $\pi_3(210)$.

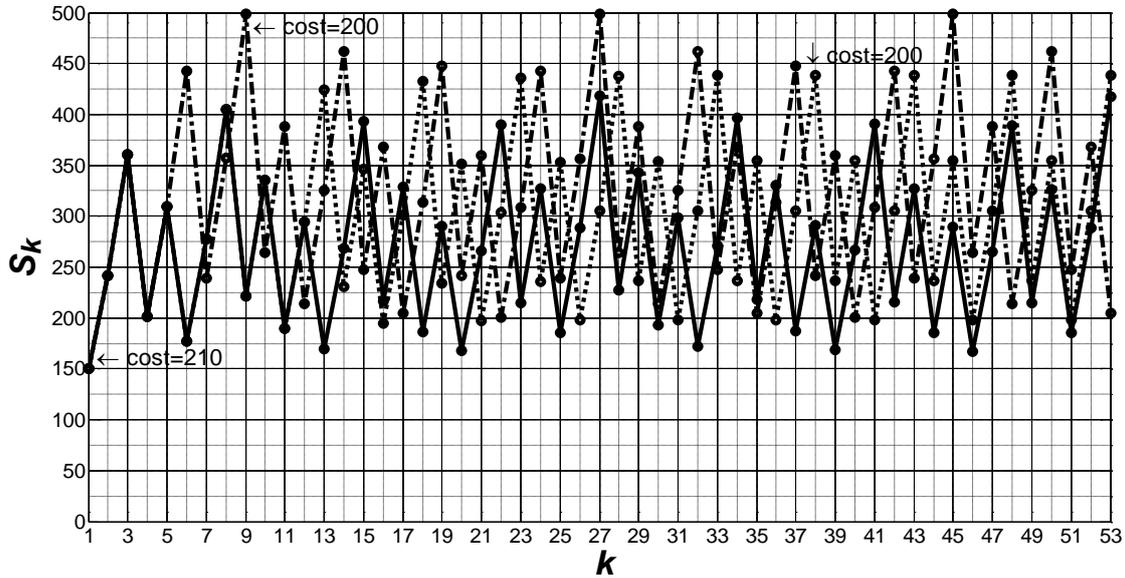

**Fig. 2** Tumor size at the beginning of each period for the optimal solution of $P_1$ (cost 200, dash-dotted line), for the optimal solution of $\pi_1(200)$ (dotted line) and for the optimal solution of $\pi_1(210)$ (solid line). In the solution of $P_1$, the maximum tumor size is equal to $S_9$=499.09, in the solution of $\pi_1(200)$ it is equal to $S_{53}$=438.42 and in the solution of $\pi_1(210)$ it is equal to $S_{27}$=418.02.

## 7 Discussion and perspective

In this article, we look at the planning of a cancer treatment over a given time horizon. The treatment consists of administering a dose of medication to the patient at certain points on the horizon. It is assumed that, in the absence of treatment, the tumor increases in size over time according to a certain growth law and, conversely, that treatment leads to a rapid reduction in tumor size according to a certain decay law without getting below a certain value. The problem considered by Bräutigam (2024) and studied in this article is to determine when treatment should be administered in order to keep tumor size below a certain threshold while minimizing the total cost of treatment.

Following Bräutigam, we approached the problem using mathematical programming, a branch of mathematics concerned with modeling and solving optimization problems. We have shown that if tumor growth follows an exponential or Gompertz law and decay consists of applying a constant reduction factor, the problem can be formulated by a mixed-integer linear program, or more precisely, by a linear program with non-negative real variables and Boolean variables. Note that this is a "direct" linear formulation obtained using the logarithmic function and not a linearization of a nonlinear problem obtained by known linearization methods, which generally increase the number of variables and/or constraints. This result enabled us to solve the problem very efficiently, using one of the many solvers available for mixed-integer linear programs. Indeed, current solvers can solve large instances of this type of program without difficulty.

The fact that we were able to efficiently solve the problem under consideration has enabled us to consider and successfully solve several interesting extensions. The first is to determine, from among all the treatment schedules that allow tumor size to be kept below a certain threshold at the lowest cost, the one that minimizes the size reached by the tumor over the horizon considered. Experiments showed that this constraint made it possible in certain cases to determine a particularly interesting schedule insofar as it enabled tumor size to be maintained at a much lower value than with the schedule determined without this constraint. This extension also makes it possible to measure precisely the marginal effect of increased treatment costs. In some cases an additional treatment dose considerably reduces the maximum tumor size over the horizon considered. Another extension of the initial problem is to consider the possibility of having several possible types of treatment. These extensions are easy to formulate from the initial problem formulation. This is one of the well-known advantages of mathematical programming: it is particularly easy to approach variants of an initial problem by adding/removing variables and/or constraints, which is not the case when the problem is approached by a specific algorithm. It should be noted, however, that "modeling" a problem using mathematical programming is not synonymous with "solving" it. Indeed, computation times can be prohibitive. The efficient resolution of the problem also means that it can be solved with a much larger number of periods. Generally speaking, the fact that we can solve the considered problem very efficiently gives us hope that we'll be able to solve problems of the same type, but which are more complicated because they are closer to real-life situations.

Being able to express the size of the tumor at the start of period $k$, $S_k$, as a function of its size at the start of period $k$-1, $S_{k-1}$, by the recurrence relations $S_k = \alpha(S_{k-1})^\beta$ if there is no treatment at the start of period $k$-1 and $S_k = \mu[\alpha(S_{k-1})^\beta]$ otherwise is the key point enabling the problem to be formulated by a linear program and therefore the key point for its efficient solution. As we have seen, the growth and decay functions considered in this article – and in (Bräutigam, 2024) – allow these expressions.

Other growth and decay functions may also need to be taken into account in a real-life context. Indeed, the effect of treatment on tumor size may be more complex than what we have considered. For example, the tumor decay caused by treatment may follow the inverse Gompertz function, which takes into account a slowdown in decay over time. Thus tumor decay would verify the equality $S_k = \mu[f(S_{k-1})]^\omega$ instead of $S_k = (1-RF)[f(S_{k-1})]$ where $\mu$ and $\omega$ are constants. It can be shown that in this case - tumor growth always following an exponential or Gompertz law - the determination of an optimal treatment schedule can be formulated by a quadratic mathematical program. To solve it, the program must either be linearized, or one of the many efficient commercial solvers that can handle mixed-integer quadratic programs directly must be used. Note also that in the case where we are interested in growth functions that cannot be expressed exactly by the equality $S_k = \alpha(S_{k-1})^\beta$ - as is the case, for example, with the Bertallanfy model (Bertallanfy, 1957) - it is possible to approach the problem by *approximating* these growth functions by such an expression.

**Declarations**

The author has no relevant financial or non-financial interests to disclose.
No funding was received to assist with the preparation of this manuscript.

**References**

Bellman, R.E. (1953). An introduction to the theory of dynamic programming.The Rand Corporation, Santa Monica, CA, Document Number R-245.

Benzekry, S., Lamont, C., Beheshti, A., Tracz, A., Ebos, J.M.L., Hlatky, L., & Hahnfeldt, P. (2014). Classical mathematical models for description and prediction of experimental tumor growth (2014). *PLOS Computational Biology,10*(8), e1003800.
https://doi.org/10.1371/journal.pcbi.1003800

Bertalanffy, L.v. (1957). Quantitative laws in metabolism and growth. *The Quarterly Review of Biology, 32*(3), 217–231. https://www.journals.uchicago.edu/doi/pdf/10.1086/401873

Bräutigam, K. (2024). Optimization of chemotherapy regimens using mathematical programming. *Computers & Industrial Engineering,191,110078*.
https://doi.org/10.1016/j.cie.2024.110078

Collins, V.P. et al. (1956). Observations on growth rates of human tumors. *Am. J. Roentgenol. Radium Ther. Nucl. Med. 76*, 988–1000

Fourer, R., Gay, D. M., & Kernighan, B. W. (1993). AMPL, a modeling language for mathematical programming. Boyd & Fraser Publishing Company.

Gerlee, P. (2013). The model muddle: in search of tumor growth laws. *Cancer Res. 73*(8), 2407-2411. https://doi.org/10.1158/0008-5472.CAN-12-4355.

Gurobi Optimization, LLC. (2023). Gurobi optimizer reference manual. Version 9.0. Gurobi Optimization, LLC. https://www.gurobi.com

Laird, A.K. (1964). Dynamics of tumor growth. *Br. J. Cancer, 13*(3), 490-502. https://doi.org/10.1038/bjc.1964.55

Mokhtari, R.B., Homayouni, T. S., Baluch, N., Morgatskaya, E., Kumar, S., Das, B., & Yeger, H. (2017). Combination therapy in combating cancer. *Oncotarget, 8*, 38022-38043. https://www.oncotarget.com/article/16723/text/

Pitchaimani, M., & Ori, G.S. (2014). Stability analysis of Gompertz tumour growth model parameters. (2014). *Journal of advances in mathematics,* 7(3), 1293-1304. https://doi.org/10.24297/jam.v7i3.7253

Sahinidis, N.V. (1996). BARON: A general purpose global optimization software package. *Journal of Global Optimization, 8*, 201-205. https://doi.org/10.1007/BF00138693

Saville, C. E., Smith, H. K., & Bijak, K. (2018). Operational research techniques applied throughout cancer care services: a review. *Health Systems*, *8*(1), 52–73. https://doi.org/10.1080/20476965.2017.1414741

Shi, J., Alagoz, O., Erenay, F.S. *et al.* (2014). A survey of optimization models on cancer chemotherapy treatment planning. *Annals of Operations Research, 221*, 331-356. https://doi.org/10.1007/s10479-011-0869-4

Yoo, T.K., Moon, H.G., Han, W., & Noh, D.Y. (2017). Time interval of neoadjuvant chemotherapy to surgery in breast cancer: how long is acceptable? *Gland Surg, 6*(1), 1-3. https://doi.org/10.21037/gs.2016.08.06

Wodarz, D., & Komarova, N.L. (2014). Dynamics of cancer: mathematical foundations of oncology. World Scientific Publishing Co. Pte. Ltd.